\documentclass[12pt, reqno,twoside]{amsart}
\usepackage[utf8]{inputenc}
\usepackage{graphicx}
\usepackage{amssymb, amsmath, amsthm, wasysym, mathrsfs}
\usepackage[in,headings]{fullpage}
\usepackage{setspace}
\usepackage{mathtools}
\usepackage{enumitem}
\usepackage{tikz}
\usepackage{comment}

\usepackage{listings}
\newtheorem{lemma}{Lemma}[section]
\newtheorem{theorem}[lemma]{Theorem}

\newtheorem{prop}[lemma]{Proposition}
\newtheorem{cor}[lemma]{Corollary}

\newtheorem{defn}[lemma]{Definition}

\theoremstyle{definition}
\newtheorem{remark}[lemma]{Remark}

\newtheorem*{lem}{Acknowledgements}

\newcommand{\Hh}{{\mathbb H}}

\newcommand{\C}{{\mathbb C}}

\newcommand{\Q}{{\mathbb Q}}

\newcommand{\Z}{{\mathbb Z}}

\newcommand{\calO}{{\mathcal O}}
\newcommand{\calP}{{\mathcal P}}
\newcommand{\calQ}{{\mathcal Q}}

\DeclareMathOperator{\rank}{rank}

\newcommand{\secref}[1]{\S\ref{#1}}

\numberwithin{equation}{section}
\numberwithin{table}{section}
\title[Some Profinitely Rigid Fibered hyperbolic 3-manifolds]{Absolute Profinite Rigidity of some closed fibered hyperbolic 3-manifolds}
\author{Tamunonye Cheetham-West}

\begin{document}

\bibliographystyle{acm}

\maketitle
\pagestyle{headings}
\begin{abstract}
    We give the first examples of closed fibered hyperbolic 3-manifolds whose fundamental groups are distinguished from every other finitely generated, residually finite group by their finite quotients. One of the examples is also the first example of a non-orientable profinitely rigid hyperbolic 3-manifold. 
\end{abstract}
\section{Introduction}
\noindent Let $G$ be a finitely generated, residually finite group and $\widehat{G}$ its profinite completion (see \secref{sec:profinite}). We say that $G$ is absolutely profinitely rigid if for any other finitely generated, residually finite group $H$ with $\widehat{H}\cong\widehat{G}$, $H\cong G$. Henceforth, an absolutely profinitely rigid group will be abbreviated to profinitely rigid, and if $G=\pi_1N$ where $N$ is a compact 3-manifold and $G$ is profinitely rigid, we will say $N$ is profinitely rigid. 
\smallbreak\noindent
 In \cite{BMRS1}, Bridson, McReynolds, Reid, and Spitler proved that the fundamental group of the Weeks manifold and some other arithmetic lattices in PSL$(2,\C)$ are profinitely rigid. The template of \cite{BMRS1} and \cite{BMRS2} to establish profinite rigidity has two main parts: first exploiting the representation rigidity of the groups concerned (using Galois rigidity, see \secref{sec:profinite}) and then the ``endgame" which uses an analysis of the subgroups of the image of a representation obtained using Galois rigidity. This endgame has now been improved and streamlined by Bridson and Reid \cite{BRPrasad} using a result of Liu \cite{Y}, that there are at most finitely many finite volume hyperbolic 3-manifolds whose fundamental groups have the same profinite completion. 
\smallbreak\noindent For the rest of the paper, $S^3_0(6_2)$ and $S^3_0(6_3)$ will denote the manifolds obtained by 0-surgery on the knots $6_2$ and $6_3$ (in the tables of Rolfsen \cite{KL}). The knots $6_2$ and $6_3$ (Figures 1 and 2) are both fibered knots of genus 2 and whose 0-surgeries are closed hyperbolic 3-manifolds that are fibered of genus 2 with infinite cyclic first homology. Furthermore, the 3-manifold $S^3_0(6_3)$ is the orientation double cover of a non-orientable 3-manifold which we will denote by $M$ for the rest of this paper. Setting $\Gamma=\pi_1(S^3_0(6_2))$ and $\Gamma'=\pi_1(S^3_0(6_3))$, the main result of this paper is 
\begin{theorem}
The groups $\Gamma$ and $\Gamma'$ are profinitely rigid. 
\end{theorem}

\noindent The proof of Theorem 1.1 follows the template of \cite{BMRS1} and \cite{BMRS2} mentioned above. Thus we first establish that $\Gamma$ is Galois rigid, then, for a finitely generated, residually finite group $\Delta$ with $\widehat{\Delta}\cong\widehat{\Gamma}$, we modify theorems from \cite{BMRS1} to produce a Zariski dense representation $\rho:\Delta\to$ PSL$(2,\C)$ and use the arithmetic of $\Gamma$ and its unique index 2 subgroup generated by squares, denoted by $\Gamma^{(2)}$ to show that $\rho(\Delta^{(2)})<\Gamma^{(2)}$ where $\Delta^{(2)}$ is the subgroup of $\Delta$ generated by squares. Using some topology and the results of \cite{BRPrasad}, one can then conclude that $\rho(\Delta^{(2)})=\Gamma^{(2)}$. It now follows that $\Delta$ is a $\Z/2\Z$ extension of $\Gamma^{(2)}$, and so the final task is to show that all $\Z/2\Z$ extensions of $\Gamma^{(2)}$ are distinguished from one another by their finite quotients, hence $\Delta\cong\Gamma$.\medbreak\noindent The strategy for showing that $\Gamma'$ is profinitely rigid is similar (see \secref{sec:sixthree} for further details). As a consequence of the classification of $\Z/2\Z$-extensions of $\Gamma'$, we are also able to prove
\begin{cor}
The non-orientable manifold $M$ is profinitely rigid. 
\end{cor}

\begin{lem}
{\it The author thanks his advisor Alan Reid for many helpful conversations and for his guidance and encouragement. The author also thanks Ryan Spitler for many helpful conversations and comments on an earlier draft of this paper}. 
\end{lem}
\begin{figure}
    \centering
    \includegraphics[scale=0.25]{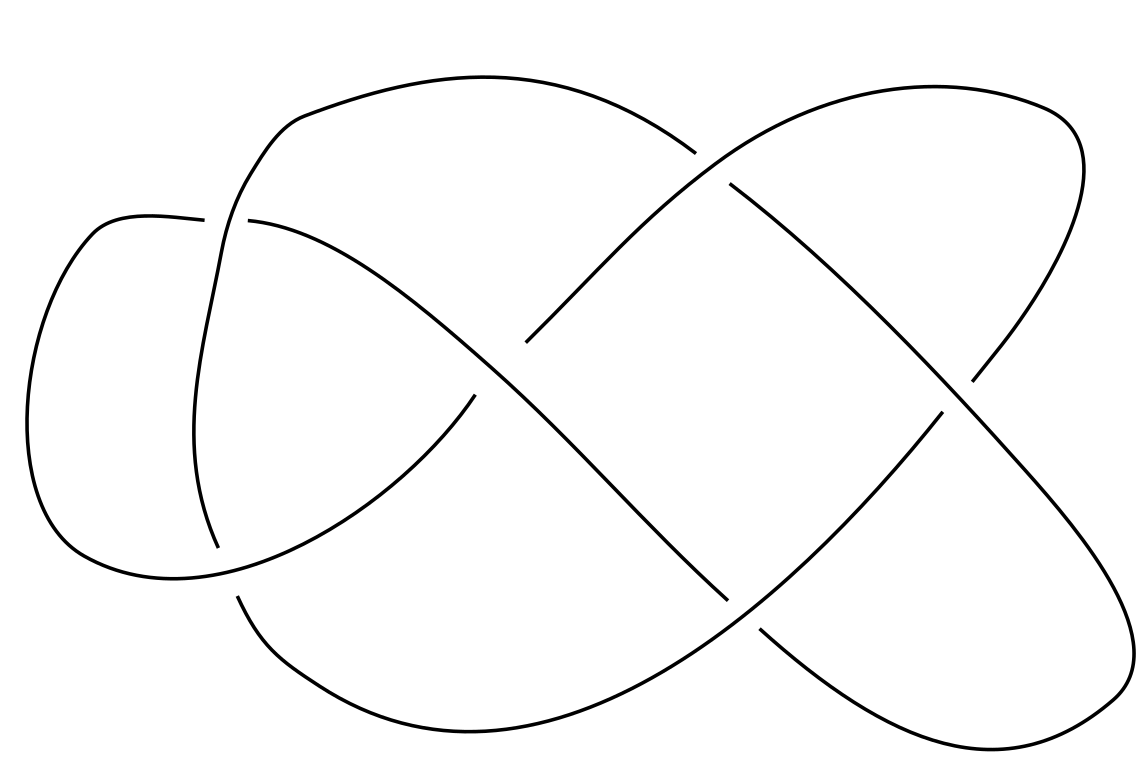}
    \caption{$6_2$}
    \label{fig:my_label2}
\end{figure}

\begin{figure}
    \centering
    \includegraphics[scale=0.25]{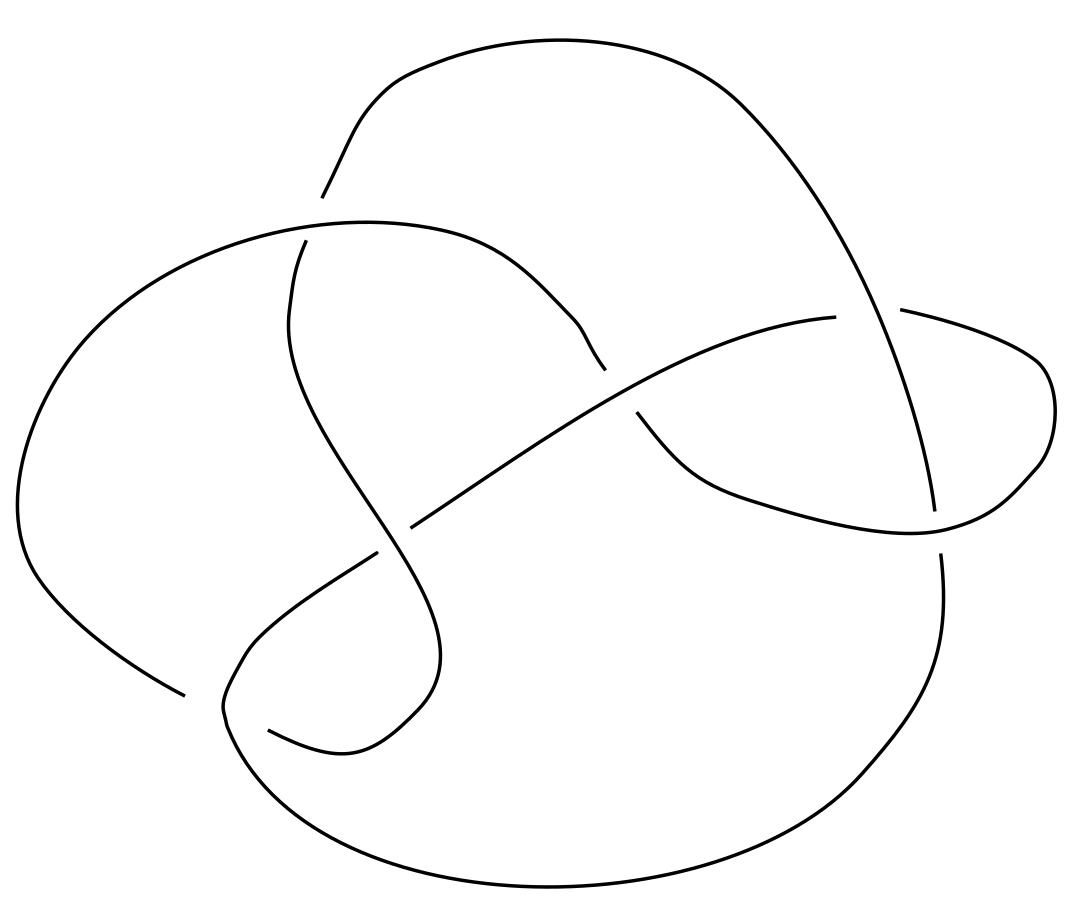}
    \caption{$6_3$}
    \label{fig:my_label}
\end{figure}

\section{Profinite Completions and Galois Rigidity}\label{sec:profinite}
\noindent The profinite completion of a group $G$ is a topological group obtained as the inverse limit of the inverse system of its finite quotients, and which is denoted by $\hat{G}$. There is a canonical map $G\to\hat{G}$ sending an element $g\in G$ to $(\Tilde{g})_{N\triangleleft_f G}$ where each $N$ is a finite-index normal subgroup of $G$ and $\Tilde{g}$ is the image of $g$ in $G/N$. The canonical map is injective if and only if $G$ is residually finite. The completion operation is functorial, taking groups to topological groups, and group homomorphisms to continuous homomorphisms of topological groups. 
\smallbreak\noindent We assume henceforth that the group $G$ is a finite covolume lattice in PSL$(2,\C)$. Let $G^{(2)}$ be the subgroup of $G$ generated by the squares of elements in $G$. This is a finite index normal subgroup with quotient the largest elementary abelian 2-group quotient of $G$. With $G$ as above, associated to $G$ are the trace-field, $k_G=\Q(tr(\gamma)\,|\,\gamma\in G)$, and the invariant trace-field, $kG=k_{G^{(2)}}$. In this setting, these fields are both number fields (\cite{MR}, Chapter 3). 
\medbreak\noindent We can also associate the following algebras to $G$
\begin{align*}
    A_0G=\lbrace\sum_{i=1}^n a_i\gamma_i\,|\,a_i\in k_G,\gamma_i\in G,n\in\Z^+\rbrace \\
AG=\lbrace\sum_{i=1}^n a_i\gamma_i\,|\, a_i\in kG, \gamma_i\in G^{(2)},n\in\Z^+\rbrace
\end{align*}
These algebras are quaternion algebras (\cite{MR}, Chapter 3) and are known as the {\it quaternion algebra} and the {\it invariant quaternion algebra} of $G$ respectively. 
\smallbreak\noindent Set $n_{k_G}=[k_G:\Q]$ which is the number of distinct embeddings of $k_G$ into the complex numbers. These embeddings give rise to pairwise non-conjugate PSL$(2,\C)$ Zariski-dense representations of $G$. Following \cite{BMRS1}, let the set of conjugacy classes of Zariski-dense PSL$(2,\C)$ representation of $G$ be $X_{zar}(G,\C)$. It follows from this discussion that $n_{k_G}\leq |X_{zar}(G,\C)|.$ 
\begin{defn}
$G$ is said to be {\it Galois rigid} whenever $n_{k_G}= |X_{zar}(G,\C)|$
\end{defn}
\noindent We can similarly define Galois rigidity for any fixed (finite covolume) representation of some abstract finitely generated group. 
\noindent For the remainder of this paper, $\Delta$ will be a finitely generated, residually finite group with $\widehat{\Delta}\cong\widehat{\Gamma}$ and $\Delta'$ will be a finitely generated, residually finite group with $\widehat{\Delta'}\cong\widehat{\Gamma'}$.
\begin{remark}
We will often use $p-$adic representations to build complex representations. To do this, we will fix some identification of $\overline{\Q_p}$ with $\C$.
\end{remark}
\subsection{Galois rigidity of $\Gamma$ and $\Gamma'$} Here we will show:
\begin{lemma}
The group $\Gamma$ is Galois rigid. 
\begin{proof}
Using Snap (\cite{Snap}) for example, it can be checked that $[k_\Gamma:\Q]=6$. Hence, to establish Galois rigidity, we need to show that $|X_{zar}(\Gamma,\C)|=6$. To that end, we shall make use of the $A-$polynomial of $6_2$ (see \cite{CCGLS}), which we denote by $A_{6_2}(L,M)$. \medbreak\noindent Recall that if a hyperbolic 3-manifold $N$ is obtained by $p/q$-surgery on a knot $J$, we can obtain a polynomial equation satisfied by the eigenvalue (denoted by $t$) of a core curve in the surgered manifold by setting $L=t^{-p},M=t^q$ in $A_{J}(L,M)$ (see \cite{LR}, for example). Since the quotient epimorphism $\pi_1(S^3\setminus 6_2)\to \Gamma$ induces an inclusion of Zariski dense characters $X_{zar}(\Gamma,\C)\hookrightarrow X_{zar}(\pi_1(S^3\setminus 6_2),\C)$ and there are no isolated points in the character variety of $6_2$ (see the Proposition on page 52 of \cite{CCGLS}), every Zariski dense PSL$(2,\C)$ representation of $\Gamma$ will be captured by the $A-$polynomial of $6_2$ by setting $L=1,M=t$.  
\smallbreak\noindent From \cite{CCGLS}, we have
\begin{multline}
    A_{6_2}(L,M)=L - 2LM^2 + L^2M^2 - M^4 + LM^4 - 3L^2M^4 + 2LM^6 + L^2M^6 - 5LM^8 + 5L^2M^8\\ - 5LM^{10}
    + 3L^2M^{10} - 3L^3M^{10} + 3LM^{12} - 12L^2M^{12} + 8L^3M^{12} - 13L^2M^{14} + 3L^3M^{14} + 3L^2M^{16}\\ - 13L^3M^{16} + 8L^2M^{18} - 12L^3M^{18} + 3L^4M^{18} - 3L^2M^{20} + 3L^3M^{20} - 5L^4M^{20} + 5L^3M^{22} - 5L^4M^{22}\\ + L^3M^{24} + 2L^4M^{24} - 3L^3M^{26} + L^4M^{26} - L^5M^{26} + L^3M^{28} - 2L^4M^{28} + L^4M^{30} 
\end{multline}
When $p=0$, and $q=1$, upon factoring $A_{6_2}(1,t)$ we obtain $$A_{6_2}(t)=(1+t^2)^5(1-3t^2+3t^4-3t^6+t^8)(1-3t^2+5t^4-5t^6+5t^8-3t^{10}+t^{12})$$ By Proposition 6.2 of \cite{CCGLS}, the factor $(1-3t^2+3t^4-3t^6+t^8)$ corresponds to reducible representations of $\pi_1(S^3\setminus 6_2)$ as well as of $\Gamma$.\medbreak\noindent We can also show by explicit computation that the factor $(1+t^2)^5$ corresponds to finite representations of the zero surgery by arguing that any representation with eigenvalue of a core curve corresponding to $t=\pm i,p=0,q=1$ will factor through the dihedral $(2,0)$ filling on the knot $6_2$. Such a representation is finite and therefore not Zariski dense. \medbreak\noindent Since $S^3_{0}(6_2)$ is a hyperbolic manifold, the eigenvalue $t_0$ of the core curve in the Zariski dense discrete faithful representation of $\Gamma$ is a zero of the last remaining factor $(1-3t^2+5t^4-5t^6+5t^8-3t^{10}+t^{12})$. Note that $t_0+t_0^{-1}\in k_\Gamma$ and has degree $6$, implying that $t_0+t_0^{-1}$ is a primitive element for $k_\Gamma$. Hence all Zariski-dense PSL$(2,\C)$ representations of $\Gamma$ (up to conjugacy) are Galois conjugates of the discrete faithful representation, and we conclude that $\Gamma$ is Galois rigid as claimed. 
\end{proof}
\end{lemma}

\begin{lemma}
The group $\Gamma'$ is Galois rigid.
\begin{proof}
As in Lemma 2.2, using Snap \cite{Snap}, we can check that $[k_{\Gamma'}:\Q]=4$. We directly verify that there are exactly four Zariski dense PSL$(2,\C)$ representations of $\Gamma'$ by directly counting the Zariski-dense PSL$(2,\C)$ characters of $\Gamma'$. From a 2-bridge presentation for $\pi_1(S^3\setminus 6_3)$, we obtain the following presentation for $\Gamma'$ $$\Gamma'=\,\langle\,a,b\,|\,wa=bw,ww^*\,\rangle\,,w=aba^{-1}b^{-1}a^{-1}bab^{-1}a^{-1}b^{-1}ab,w^*=bab^{-1}a^{-1}b^{-1}aba^{-1}b^{-1}a^{-1}ba$$ where $ww^*$ represents a homological longitude. Any Zariski dense $PSL(2,\C)$ character of $\Gamma'$ corresponds to a Zariski dense PSL$(2,\C)$ character of $\pi_1(S^3\setminus 6_3)$ arising from a representation that maps the word $ww^*$ to the identity in PSL$(2,\C)$. Furthermore, for any such Zariski dense representation $f:\Gamma'\to $PSL$(2,\C)$, we can conjugate such that $f(a)=\left(\begin{smallmatrix}p&1\\0& 1/p\end{smallmatrix}\right)$, $f(b)=\left(\begin{smallmatrix}p&0\\r& 1/p\end{smallmatrix}\right)$. With this setup, one can verify that there are four Zariski dense representations of $\Gamma'$. A Mathematica \cite{Mathematica} notebook where this computation is done is available upon request from the author. 
\end{proof}
\end{lemma}

\section{Constructing PSL(2,$\C$) Representations For $\Delta$}
\subsection{The arithmetic structure of $\Gamma$}\label{sec:arith}
The group $\Gamma$ is arithmetic and lies in the commensurability class of the Weeks manifold (see \cite{BMRS1}, page 709) which has (invariant) trace field $k=\Q(\theta)$ with $\theta^3-\theta^2+1=0$ (with ring of integers $R_k$), a field of discriminant -23 and invariant quaternion algebra $B/k$ that is ramified at the unique real place and the unique place of norm 5 in $k$. The algebra $B/k$ is {\it locally uniform} in the sense of \cite{BMRS1} which means that for any two places $\mu,\nu$ of $k$ for which $k_\mu\cong k_\nu$, $B_\mu\cong B_\nu$.Moreover, the group $\Gamma$ has integral traces and is not derived from a quaternion algebra. However $\Gamma^{(2)}<\Gamma^1_{\calO}$, the subgroup of units of the unique (conjugacy class of) maximal order $\calO$ of $B$. A fact that is helpful for the Magma computation in Section 7.2 is that $\Gamma^1_\calO$ is the orbifold fundamental group of the $(3,0)$ Dehn filling on the knot $5_2$ with finite presentation (\cite{SnapPy}) $$\Gamma^1_\calO=\,\langle\,a,b\,|\,a^3,ab^{-1}ab^2a^{-1}b^{-2}a^{-1}b^2\,\rangle$$ 
\smallbreak\noindent Following Section 9.1 of \cite{BMRS1}, we observe that there are two primes, denoted by $\calP$ and $\calQ$, both of norm 23, one ramified and one split (thus $k_\calP$ and $k_\calQ$ are non-isomorphic extensions of $\Q_{23}$) and the reduction homomorphisms $R_k\to R_k/\calP\cong\Z/23\Z$ and $R_k\to R_k/\calQ\cong \Z/23\Z$ induce two reduction epimorphisms $\Gamma^1_\calO\to $PSL$(2,\Z/23\Z)$. Furthermore, we can check, using Magma \cite{Magma} (see Section 7.1) that these are the only two PSL$(2,\Z/23\Z)$ quotients of $\Gamma^1_\calO$. We will also make use of a fact from Section 9.1 of \cite{BMRS1} that there is a unique conjugacy class of (maximal) index 24 subgroup in PSL$(2,\Z/23\Z)$ as follows:\smallbreak\noindent For one of these reduction homomorphisms corresponding to $\calP$, say, which we will denote by $\phi_\calP:\Gamma^1_\calO\to $PSL$(2,\Z/23\Z)$, the image of $\Gamma^{(2)}$ has index 24 in PSL$(2,\Z/23\Z)$, and the core of $\Gamma^{(2)}$ in $\Gamma^1_\calO$ is exactly $\ker\phi$ (see Section 7.2). It follows that up to conjugacy, $\Gamma^{(2)}$ is the full preimage of a fixed (maximal) index 24 subgroup of PSL$(2,\Z/23\Z)$ under $\phi_\calP$. Observing that $\phi_{\calP}$ factors through PSL$(2,R_{k_\calP})$, where $R_{k_\calP}$ is the non-Archimedean completion of $R_k$ at the place corresponding to $\calP$ and the ring of integers of $k_\calP$ the $\calP$-adic completion of $k$ (see the remarks following Theorem 0.6.6 of \cite{MR}), we have that there is a corresponding maximal open subgroup of index 24 in PSL$(2,R_{k_\calP})$ denoted by $H$ (unique up to conjugation in PSL$(2,R_{k_\calP})$) whose preimage in $\Gamma^1_\calO$ is $\Gamma^{(2)}$. 
 
\subsection{Running the BMRS program on $\Gamma$}
We will apply Theorem 4.8 of \cite{BMRS1}. This theorem implies that since $\Gamma$ is Galois rigid by Lemma 2.2, and with $\Delta$ as stated, we obtain a Zariski dense representation $\rho:\Delta\to $PSL$(2,\C)$ such that $k_{\rho(\Delta)}$ and $k_\Gamma$ are arithmetically equivalent fields. A consequence of Theorem 4.8 in \cite{BMRS1} that is applicable for our setting is:
\begin{lemma}
With $\Gamma,\Delta$ as above, there is a representation $\rho:\Delta\to$ PSL$(2,\C)$ with Zariski-dense image such that 
\begin{enumerate}
    \item $\rho(\Delta)$ is Galois rigid;
    \item For $V_{k_{\rho(\Delta)}}$ and $V_{k_\Gamma}$ the sets of finite places of $k_{\rho(\Delta)}$ and $k_\Gamma$ respectively, there is a bijection $\tau: V_{k_\Gamma}\to V_{k_{\rho(\Delta)}}$ such that the local fields $k_{\rho(\Delta)_{\tau(\nu)}}$ and $k_{\Gamma_\nu}$ are isomorphic for every $\nu\in V_{k_{\Gamma}}$. 
\end{enumerate}
\end{lemma}
\noindent Since $\Gamma$ is not derived from a quaternion algebra, it will be convenient to pass to the groups $\Gamma^{(2)}$ and $\Delta^{(2)}$. To this end, the following lemma compares $k$ and $k_{\rho(\Delta^{(2)})}$.
\begin{lemma}
The image, $\rho(\Delta^{(2)})$, of the restriction of $\rho$ to the unique index 2 subgroup $\Delta^{(2)}<\Delta$ has trace field $k$. 
\begin{proof}
Let $l=k_{\rho(\Delta^{(2)})}$. For every rational prime $p$, we will build a bijection $\tau^p_2$ from the finite places of $k$ lying above $p$ denoted $V^p_k$ to the finite places of $l$ lying above $p$ denoted $V^p_l$. For an element $\nu\in V^p_k$, pick a prime $\mu$ in $V^p_{k_\Gamma}$ that lies above $\nu$. Consider the prime $\mu'\in V^p_{k_{\rho(\Delta)}}$ with $\tau(\mu)=\mu'$ using the bijection from Lemma 3.1 (2). The prime $\mu'$ lies above a prime $\nu'$ of $l$. Set $\tau^p_2(\nu)=\nu'$. 
\medbreak\noindent We show that $k_\nu\cong l_{\nu'}$ for all finite places $\nu$ of $k$. 
Since $k_{\Gamma,\mu}\cong k_{\rho(\Delta),\mu'}$ by Lemma 3.1 (2), we have that $k_\nu\cong k_{\Gamma^{(2)},\mu}\cong k_{\rho(\Delta)^{(2)},\mu'}\cong l_{\nu'}$. The middle isomorphism $k_{\Gamma^{(2)},\mu}\cong k_{\rho(\Delta)^{(2)},\mu'}$ follows from the way $\rho$ is constructed in \cite{BMRS1} as a map obtained by composing the canonical injection $\Delta\hookrightarrow\widehat{\Gamma}$ with a bounded, continuous PSL$(2,\overline{\Q_p})-$representation of $\widehat{\Gamma}$. From this construction, it follows that the traces of $\Gamma^{(2)}$ and $\rho(\Delta^{(2)})$ will have the same topological closure at every corresponding pair of places $\mu,\tau(\mu)$ of $k_\Gamma$ and $k_{\rho(\Delta)}$ respectively. 
\smallbreak\noindent To see that $\tau_2^p$ is a well-defined bijection, if we use some other prime $\mu''$ lying above $\nu$, then, the non-Archimedean completions at the places $\mu$ and $\mu''$ have an isomorphic subfield corresponding to the completion (of $k$) at $\nu$. The field $k_\nu$ is exactly the topological closure of $k$ in $k_{\Gamma,\mu}$ under the completion field homomorphism. For $\mu''$ also lying above $\nu$, $k_{\nu'}$ is also exactly the topological closure of $k$ in $k_{\Gamma,\mu''}$.
\medbreak\noindent Having constructed bijections $\tau_2^p$ for every finite rational prime $p$, we can form a bijection between the finite places by setting $\tau_2=\sqcup_p \tau_2^p:V_k\to V_l$ and noting that $V_l=\sqcup_p V^p_l$ and $V_k=\sqcup_p V^p_k$. From Lemma 4.9 of \cite{BMRS1}, it follows that $k_{\rho(\Delta^{(2)})}$ and $k$ are arithmetically equivalent. Since arithmetically equivalent fields of small degree ($\leq 7$) are isomorphic \cite{P}, we conclude that $k_{\rho(\Delta^{(2)})}\cong k$.
\end{proof}
\end{lemma}
\begin{prop}
The algebra $A_0\rho(\Delta^{(2)})\cong B$.
\begin{proof}
By construction, $B_\nu\cong A_0\rho(\Delta^{(2)})_{\tau_2(\nu)}$ at every finite place $\nu$ of $k$. Since $B$ is ramified at the unique finite place of norm 5 and since $B$ and $A_0\rho(\Delta^{(2)})$ are quaternion algebras over the same cubic field, they are isomorphic by Corollary 4.11 (III) of \cite{BMRS1} because $B$ is locally uniform (as previously noted in \secref{sec:arith}). 
\end{proof}
\end{prop}
\noindent We can now proceed to strengthen the conclusions of Lemma 3.1 and Lemma 3.2 as follows 
\begin{lemma}
With $\Delta$ as above, there is a representation $\rho:\Delta\to$ PSL$(2,\C)$ such that $\rho(\Delta^{(2)})<\Gamma^{(2)}$.
\begin{proof}
Recall the prime $\calP$ from \secref{sec:arith}. Choose a prime $\calP'$ of $k_\Gamma$ lying above $\calP$. Since $\Gamma$ has integral traces, we obtain a bounded PSL$(2,\overline{\Q_{23}})$ representation $\psi:\Gamma\to $PSL$(2,\overline{\Q_{23}})$ by completing $k_\Gamma$ at $\calP'$. The closure $\overline{\psi(\Gamma)}$ is profinite and so the homomorphism $\psi$ factors through the canonical homomorphism $\Gamma\hookrightarrow\widehat{\Gamma}$. Since $\widehat{\Delta}\cong\widehat{\Gamma}$, we can now compose with the canonical homomorphism $\Delta\hookrightarrow\widehat{\Delta}$ to obtain a Zariski dense homomorphism $\Tilde{\rho}:\Delta\to $PSL$(2,\overline{\Q_{23}})$ with $\overline{\Tilde{\rho}(\Delta)}=\overline{\psi(\Gamma)}$. Thus, we obtain a homomorphism $\rho:\Delta\to $PSL$(2,\C)$. \medbreak\noindent For every prime $\calQ$ of $k$ and for every prime $\calQ'$ of $k_{\rho(\Delta)}$ lying above $\calQ$, we get a bounded representation of $\Delta$ into PSL$(2,R_{k_{\Gamma},{\calQ'}})$ where $R_{k_{\Gamma},{\calQ'}}$ is the local ring of $k_{\Gamma_{\calQ'}}$, and so $\rho(\Delta^{(2)})$ has integral traces because $$R_{k_\Gamma}=\bigcap_{\calQ'\,\text{a finite place of }k_\Gamma}\enspace(k_\Gamma\cap R_{k_{\Gamma},{\calQ'}})$$ Applying Lemma 3.2, Proposition 3.3, and Corollary 4.11 (III) of \cite{BMRS1}, it follows that $\rho(\Delta^{(2)})<B$. Since $\rho$ is an integral representation and $B$ has type number 1 (having a unique conjugacy class of maximal order), we can conjugate the representation such that $\rho(\Delta^{(2)})<\Gamma^1_\calO$ as claimed.
\medbreak\noindent We now use the information in the latter part of \secref{sec:arith} as follows: by the construction of $\rho$ using the prime $\calP$, $\rho(\Delta^{(2)})$ has the same topological closure with $\Gamma^{(2)}$ at their corresponding places (using the bijection $\tau_2$ from the proof of Lemma 3.2). Thus, at the place $\tau_2(\calP)$, $\rho(\Delta^{(2)})<\Gamma^1_\calO$ has the same bounded topological closure as $\Gamma^{(2)}$. Since (as noted in \secref{sec:arith}) $\Gamma^{(2)}$ is exactly the preimage of a representative of the unique conjugacy class of index 24 open subgroups of PSL$(2,R_{k_\calP})$, where $k_\calP$ is the non-archimedean completion of the field at that place, we have that $\rho(\Delta^{(2)})<\Gamma^{(2)}$.
\end{proof}
\end{lemma}
\smallbreak\noindent 
\section{Constructing PSL(2,$\C$) representations for $\Delta'$}
\subsection{Arithmetic information for $\Gamma'$} Using Snap \cite{Snap}, it can be checked that the group $\Gamma'$ is arithmetic, has invariant trace field $\Q(\sqrt{-3})$ and invariant quaternion algebra $A\Gamma'/\Q(\sqrt{-3})$
which is ramified at the unique places of norms 3 and $5^2$. It follows that $A\Gamma'$ is locally uniform. While $\Gamma'$ is not derived from a quaternion algebra, $\Gamma'^{(2)}$ is. In particular, from Snap \cite{Snap} we can compute $vol(\Hh^3/\Gamma'^{(2)})$ ($\Hh^3/\Gamma'^{(2)}$ is the unique index 2 cover of the 0-surgery on $6_3$) and observe that it is equal to the covolume of the unit group of a maximal order of $A\Gamma'$ computed using Theorem 11.1.3 of \cite{MR}. Thus, $\Gamma'^{(2)}$ is exactly the unit group $\Gamma^1_{\calO'}$ of a maximal order $\calO'$ of the invariant quaternion algebra $A\Gamma'$. We also note that the quaternion algebra $A\Gamma'$ has type number 1 and so there is a unique conjugacy class of maximal order in $A\Gamma'$.    
\subsection{Producing a Zariski dense representation of $\Gamma'$} Since $\Gamma'$ is Galois rigid, we have the following 

\begin{lemma}
There is a representation $\rho':\Delta'\to $PSL$(2,\C)$ such that $\rho'(\Delta'^{(2)})<\Gamma'^{(2)}$. 
\begin{proof}
Using analogues of Lemmas 2.3, 3.1, and 3.2 to $\Gamma',\Delta'$, it follows that we can choose $\rho'$ such that the trace field of $\rho'(\Delta')$ is $k_{\Gamma'}$ and the trace field of $\rho'(\Delta'^{(2)})$ is $k\Gamma'$. In a similar manner to the proof of Lemma 3.4, we show that $\rho'$ is an integral representation. Also, $A\rho'(\Delta'^{(2)})$ is locally equivalent to $A\Gamma'$ as quaternion algebras over the same field by an analog of Proposition 3.3. Therefore, since $A\Gamma'$ is locally uniform, we apply Corollary 4.11 (II) of \cite{BMRS1} to show that $A\rho'(\Delta'^{(2)})\cong A\Gamma'$. Since the restriction of $\rho'$ to $\Delta'^{(2)}$ is an integral representation into a type number 1 quaternion algebra and $\Gamma'^{(2)}=\Gamma^1_{\calO'}$, we can can conjugate the representation such that $\rho'(\Delta'^{(2)})<\Gamma'^{(2)}$ as claimed. 
\end{proof}
\end{lemma}

\section{Profinite Rigidity of $\pi_1(S_0^3(6_2))$}\label{sec:sixtwo}
\noindent To prove that $\rho(\Delta^{(2)})=\Gamma^{(2)}$, we use Proposition 5.1 from \cite{BRPrasad} stated here for convenience. We also use the fact that since $H_1(\Gamma,\Z)\cong\Z$ and $H_1(\Gamma^{(2)},\Z)\cong\Z\oplus \Z/11\Z$, the same holds for $H_1(\Delta,\Z)$ and $H_1(\Delta^{(2)},\Z)$ respectively. 

\begin{prop}
Let $G$ be a Kleinian group of finite covolume. If $\Lambda$ is a finitely generated, residually finite group with a homomorphism $\Lambda\to G$ whose image is a proper subgroup of finite index then $\widehat{G}$ and $\widehat{\Lambda}$ are not isomorphic.
\end{prop}
\begin{lemma}
For a group $\Delta$ as above, with $\rho:\Delta\to $PSL$(2,\C)$ as constructed in Lemma 3.3 above, $\rho(\Delta^{(2)})=\Gamma^{(2)}$ and $\rho|_{\Delta^{(2)}}$ is an isomorphism. 
\begin{proof}
If $\rho(\Delta^{(2)})<\Gamma^{(2)}$ is an infinite-index subgroup, the manifold $\mathbb{H}^3/\rho(\Delta^{(2)})$ has a compact core with some boundary component having genus $\geq 2$. By Poincar{\'e}-Lefschetz duality, we deduce that $b_1(\rho(\Delta^{(2)}))=\rank(H_1(\rho(\Delta^{(2)}),\Q))>1$, a contradiction since $b_1(\Delta^{(2)})=1$. Thus $\rho(\Delta^{(2)})<\Gamma^{(2)}$ is a finite index subgroup and we can apply Proposition 5.1 above to conclude that $\rho(\Delta^{(2)})=\Gamma^{(2)}$. Since $\widehat{\Gamma^{(2)}}\cong\widehat{\Delta^{(2)}}$, it follows from the Hopfian property of finitely generated profinite groups (Proposition 2.5.2 of \cite{RZ}) that $\rho$ is injective, hence $\Delta^{(2)}\cong\Gamma^{(2)}$. 
\end{proof}

\end{lemma}

\noindent So far, we have established that for a group $\Delta$ with $\widehat{\Delta}\cong\widehat{\Gamma}$, $\Delta^{(2)}\cong\Gamma^{(2)}$. Thus, $\Delta$ is a $\Z/2\Z$ extension of $\Gamma^{(2)}$. To argue that $\Delta\cong\Gamma$, it is sufficient to show that all $\Z/2\Z$ extensions of $\Gamma^{(2)}$ are distinguished from one another by their finite quotients. 
\begin{lemma}
Let $\Delta$ be as above. Then $\Delta\cong\Gamma$.
\begin{proof}
The group $\Gamma$ is {\it good} in the sense of Serre (see Chapter 2, Exercise 2 of \cite{Serre} for a definition) and so $\widehat{\Gamma}$ is torsion-free (see Section 2.3 and Lemma 2.5 of \cite{BRPrasad} for example). Since $\widehat{\Delta}\cong\widehat{\Gamma}$, it follows that $\Delta$ is torsion-free . Furthermore, since $\Delta$ is torsion free, we can apply Proposition 4.1 of \cite{BRPrasad} to show that $\Delta$ is a lattice in PSL$(2,\C)$ as well as a $\Z/2\Z$ extension of $\Gamma^{(2)}$ with $H_1(\Delta,\Z)\cong\Z$ as noted earlier.   
\medbreak\noindent
Using SnapPy \cite{SnapPy} and as a consequence of Mostow rigidity, we check that $Out(\Gamma^{(2)})$ is $(\Z/2\Z)^3$. It follows that $\Gamma^{(2)}$ is index 8 in $N(\Gamma^{(2)})$, its normalizer in $Isom(\Hh^3)$. The group $\Gamma^{(2)}$ is the unique index 2 subgroup of $\Gamma$ since $H_1(\Gamma,\Z)\cong\Z$, and it follows that the normalizer $N(\Gamma)$ of $\Gamma$ in $Isom(\Hh^3)$ is a subgroup of $N(\Gamma^{(2)})$. Another SnapPy calculation shows that $Out(\Gamma)$ is $(\Z/2\Z)^2$, and so $N(\Gamma)=N(\Gamma^{(2)})$, since $\Gamma$ is an index 4 subgroup of its normalizer in $Isom(\Hh^3)$. 
\medbreak\noindent We now find a finite presentation for $N(\Gamma)=N(\Gamma^{(2)})$ with which we compute the integral first homology group of all non-trivial extensions of $\Gamma^{(2)}$ and show that there is a unique extension with infinite cyclic first homology and that is $\Gamma$. To find this finite presentation, we use the fact that all the symmetries of $\Gamma$ arise from orientation preserving symmetries of the 2-bridge knot $6_2$ (see Theorem 4.1 of \cite{S}, for example), and these extend to symmetries of $S^3_0(6_2)$ for homological reasons. In particular, the $(\Z/2\Z)^2$ group comes from two orientation-preserving order 2 symmetries of $6_2$ which we will denote $\sigma$ and $\tau$. A 2-bridge presentation for $\pi_1(S^3\setminus 6_2)$ is given by
$$\pi_1(S^3\setminus 6_2)=\,\langle\, a,b\,|\,wa=bw\,\rangle\,, w=abab^{-1}a^{-1}b^{-1}a^{-1}bab$$
we can assume that $\sigma$ sends $a\to b$, $b\to a$ and $\tau$ sends $a\to a^{-1}$ and $b\to b^{-1}$. From this, we can get a finite presentation for $\Gamma$ (the word $a^{-4}ww^*$ represents a homological longitude for $6_2$)
$$\Gamma=\,\langle\,a,b\,|\,wa=bw,a^{-4}ww^*\,\rangle\,,w^*=baba^{-1}b^{-1}a^{-1}b^{-1}aba$$
and then we can write a presentation for $N(\Gamma)$ as follows; $$\langle\,a,b,\sigma,\tau\,|\,wa=bw,a^{-4}ww^*,\sigma\tau\sigma^{-1}\tau^{-1},\sigma a\sigma^{-1}b^{-1},\sigma b\sigma^{-1}a^{-1},\tau a\tau^{-1}a,\tau b\tau^{-1}b,\sigma^2,\tau^2,(\sigma\tau)^2\,\rangle$$
With this presentation of $N(\Gamma)$, we use a Magma script (in Section 7.3) to compute the first integral homology of the $\Z/2\Z$ extensions of $\Gamma^{(2)}$ and $\Gamma$ is the only extension of $\Gamma^{(2)}$ with integral first homology $\Z$, completing the proof of Lemma 5.3 and the first part of Theorem 1.1. 
\end{proof}
\end{lemma}

\section{Profinite Rigidity of $\pi_1(S_0^3(6_3))$ and $\pi_1M$}\label{sec:sixthree}
\noindent In this section, we will use the fact that $H_1(\Gamma',\Z)\cong\Z$ and therefore $H_1(\Delta',\Z)\cong\Z$. Similarly, since $H_1(\Gamma'^{(2)},\Z)\cong\Z\oplus\Z/13\Z$, we have $H_1(\Delta'^{(2)},\Z)\cong\Z\oplus\Z/13\Z$. We follow the proof in \secref{sec:sixtwo}. First, the proof of Lemma 5.2 applies to show that $\rho'(\Delta'^{(2)})\cong\Gamma'^{(2)}$. 
\noindent We can now conclude the proof of Theorem 1.1 with the following 
\begin{lemma}
The group $\Delta'\cong\Gamma'$. 
\begin{proof}
The proof proceeds exactly as before to first show that $\Gamma'^{(2)}$ has index 8 in $N(\Gamma'^{(2)})$, its normalizer in PSL$(2,\C)$, and that $\Gamma'$ and $\Gamma'^{(2)}$ have the same normalizer in PSL$(2,\C)$. 
\medbreak\noindent As before, we find a finite presentation for $N(\Gamma')=N(\Gamma'^{(2)})$ with which we compute the integral first homology group of all non-trivial extensions of $\Gamma'^{(2)}$ and show that there is a unique extension with infinite cyclic first homology and that is $\Gamma'$. To find this finite presentation, we use the fact that all the orientation preserving symmetries of $\Gamma'$ arise from orientation preserving symmetries of the 2-bridge knot $6_3$. In particular, the orientation preserving $(\Z/2\Z)^2$ group comes from two order 2 symmetries of $6_3$ which we will denote $\sigma$ and $\tau$. A 2-bridge presentation for $\pi_1(S^3\setminus 6_3)$ is given by
$$\pi_1(S^3\setminus 6_3)=\,\langle\, a,b\,|\,wa=bw\,\rangle\,, w=aba^{-1}b^{-1}a^{-1}bab^{-1}a^{-1}b^{-1}ab$$
we can assume that $\sigma$ sends $a\to b$, $b\to a$ and $\tau$ sends $a\to a^{-1}$ and $b\to b^{-1}$. From this, we can get a finite presentation for $\Gamma'$ (the word $ww^*$ represents a homological longitude for $6_3$)
$$\Gamma'=\,\langle\,a,b\,|\,wa=bw,ww^*\,\rangle\,,w^*=bab^{-1}a^{-1}b^{-1}aba^{-1}b^{-1}a^{-1}ba$$
and then we can write a presentation for $N(\Gamma')$; $$\langle\,a,b,\sigma,\tau\,|\,wa=bw,ww^*,\sigma\tau\sigma^{-1}\tau^{-1},\sigma a\sigma^{-1}b^{-1},\sigma b\sigma^{-1}a^{-1},\tau a\tau^{-1}a,\tau b\tau^{-1}b,\sigma^2,\tau^2,(\sigma\tau)^2\,\rangle$$
With this presentation of $N(\Gamma')$, we use a Magma script (in Section 7.4) to compute the first integral homology of the $\Z/2\Z$ extensions of $\Gamma'^{(2)}$ and $\Gamma'$ is the only extension of $\Gamma'^{(2)}$ with integral first homology $\Z$, concluding the proof of this lemma and Theorem 1.1.
\end{proof}
\end{lemma}
\subsection{A profinitely rigid non-orientable 3-manifold}
The aim of this subsection is to deduce Corollary 1.2. The group $\pi_1M$ is a $\Z/2\Z$-extension of $\Gamma'$ and so it is sufficient to distinguish the $\Z/2\Z$-extensions of $\Gamma'$ by integral homology i.e. to prove that $\pi_1M$ is the only extension with integral first homology $\Z$. Let $G$ be a finitely generated, residually finite group with $\widehat{G}\cong\widehat{\pi_1M}$. The following lemma completes the proof of Corollary 1.2
\begin{lemma}
The group $G$ as above is a $\Z/2\Z$-extension of $\Gamma'$ with $H_1(G,\Z)\cong\Z$. Furthermore, $\pi_1M$ is the unique $\Z/2\Z$-extension of $\Gamma'$ with first integral homology $\Z$. Thus, $G\cong\pi_1M$.
\begin{proof}
The profinite group $\widehat{G}$ has a unique index two subgroup isomorphic to $\widehat{\Gamma'}$. By Theorem 1.1, $G$ has a unique index two subgroup isomorphic to $\Gamma'$, and so $G$ is a $\Z/2\Z$ extension of $\Gamma'$. Since $H_1(M,\Z)\cong\Z$, the first integral homology of $G$ is $\Z$ as well. 
\medbreak\noindent The 0-surgery homomorphism from $\pi_1(S^3\setminus 6_3)$ to $\Gamma'$ induces an isomorphism between $Out(\pi_1(S^3\setminus 6_3))$ and $Out(\Gamma')$. Thus, every $\Z/2\Z$-extension of $\Gamma'$ corresponds to a Dehn filling on a $\Z/2\Z$-extension of $\pi_1(S^3\setminus 6_3)$. Using the following presentation of the maximal discrete subgroup $N(\pi_1(S^3\setminus 6_3))$ containing $\pi_1(S^3\setminus 6_3)$ in $Isom(\Hh^3)$ obtained from Page 111 of \cite{riley82} $$N(\pi_1(S^3\setminus 6_3))=\,\langle a,b,c,d,e\,|\,bdba,c^{-1}dbe^{-1}aea,a^2,b^2,(bc)^2,c^2,d^2,ed^{-1}ce\,\rangle$$ one can check (see Section 7.5) that there are only two index four subgroups of $N(\pi_1(S^3\setminus 6_3))$ with positive first Betti number. We can check that one of them is orientable with integral first homology $\Z\oplus\Z/2\Z$ and the other is non-orientable with first homology $\Z$. The latter manifold is m018 in the non-orientable census. Thus, if $G$ is not the fundamental group of a Dehn filling on m018 then $G$ has integral homology  $\Z\oplus\Z/2\Z$, contradicting $H_1(G,\Z)\cong H_1(\pi_1M,\Z)\cong\Z$. Alternatively, one can argue that if $G\not\cong\pi_1M$, then $G$ is the fundamental group of a closed orientable 3-manifold. Since $G$ and $\pi_1M$ are cohomologically good, $H^3(G,\Z/3\Z)\cong H^3(\widehat{\pi_1M},\Z/3\Z)\cong H^3(\pi_1M,\Z/3\Z)$. However, $H^3(G,\Z/3\Z)\cong\Z/3\Z$, $H^3(\pi_1M,\Z/3\Z)=0$ and we obtain a contradiction, concluding the proof of the corollary. 
\end{proof}
\end{lemma}

\section{Routines used in the paper}
\subsection{Magma script for checking PSL$(2,\Z/23\Z)$ quotients for $\Gamma^1_\calO$}
\begin{verbatim}
> M<a,b>:=FPGroup<a,b|a^3,a*b^-1*a*b^2*a^-1*b^-2*a^-1*b^2>; 
> L2Quotients(M);
[
    L_2(infty^3)
]
> quot:=L2Quotients(M);
> q:=quot[1];
> SpecifyCharacteristic(q,23);
[
    L_2(23),
    L_2(23)
]
>
\end{verbatim}

\subsection{Verifying the congruence information for $\Gamma^{(2)}$}
\begin{verbatim}
> g:=PSL(2,23);
> homs:=Homomorphisms(M, g : Surjective);
> #homs;
4
> list:=LowIndexSubgroups(M,<24,24>);
> for j in [1 .. #list] do
for> print j, AQInvariants(list[j]);
for> end for;
1 [ 30 ]
2 [ 11, 0 ]
3 [ 5, 30 ]
4 [ 3, 6, 6 ]
5 [ 66 ]
6 [ 3, 6, 6 ]
7 [ 7, 42 ]
8 [ 3, 6, 0 ]
9 [ 2, 2, 6 ]
10 [ 3, 6, 0 ]
11 [ 2, 2, 6 ]
> sub:=list[2];          
> Index(g,homs[1](sub));
1
> Index(g,homs[2](sub));
1
> Index(g,homs[3](sub));
24
> Index(g,homs[4](sub));
24
>     
\end{verbatim}

\subsection{Magma script for distinguishing non-trivial $\Z/2\Z$ extensions of $\Gamma^{(2)}$}
\begin{verbatim}
> N<a,b,s,t>:=FPGroup<a,b,s,t|s^2,t^2,(s*t)^2,s*a*s^-1*b^-1,s*b*s^-1*a
^-1,t*a*t^-1*a,t*b*t^-1*b,a^-4*(a*b*a*b^-1*a^-1*b^-1*a^-1*b*a*b)*(b*a*
b*a^-1*b^-1*a^-1*b^-1*a*b*a),(a*b*a*b^-1*a^-1*b^-1*a^-1*b*a*b)*a*(b^-1
*a^-1*b^-1*a*b*a*b*a^-1*b^-1*a^-1)*b^-1>;
> list:=LowIndexSubgroups(N,<4,4>);
> for j in [1 .. #list] do
for> if IsNormal(N,list[j]) then
for|if> print AQInvariants(list[j]);
for|if> end if;
for> end for;
[ 2, 22 ]
[ 2, 2 ]
[ 2, 2 ]
[ 2, 22 ]
[ 11, 0 ]
[ 2, 0 ]
[ 0 ]
> 

\end{verbatim}

\subsection{Magma script for distinguishing $\Z/2\Z$ extensions of $\Gamma'^{(2)}$ that are lattices in PSL(2,$\C$).}
\begin{verbatim}
> N<a,b,s,t>:=FPGroup<a,b,s,t|b*(a*b*a^-1*b^-1*a^-1*b*a*b^-1*a^-1*b^-1*a*b)
*a^-1*(b^-1*a^-1*b*a*b*a^-1*b^-1*a*b*a*b^-1*a^-1),
(a*b*a^-1*b^-1*a^-1*b*a*b^-1*a^-1*b^-1*a*b)*(b*a*b^-1*a^-1*b^-1*a*b*a^-1*b
^-1*a^-1*b*a),s*a*s^-1*b^-1,s*b*s^-1*a^-1,t*a*t^-1*a,t*b*t^-1*b,s^2,
t^2,(s*t)^2>;
> AbelianQuotientInvariants(N);
[ 2, 2, 2 ]
> k:=LowIndexSubgroups(N,<4,4>);
> for j in [1 .. #k] do
for> if IsNormal(N,k[j]) then
for|if> print AbelianQuotientInvariants(k[j]);
for|if> end if;
for> end for;
[ 2, 26 ]
[ 2, 2 ]
[ 2, 2 ]
[ 2, 26 ]
[ 13, 0 ]
[ 2, 0 ]
[ 0 ]
\end{verbatim}

\subsection{Magma script for calculating abelianizations of $\Z/2\Z$ extensions of $\pi_1(S^3\setminus 6_3)$.}
\begin{verbatim}
> S<a,b,c,d,e>:=FPGroup<a,b,c,d,e|b*d*b*a,c^-1*d*b*e^-1*a*e*a,
a^2,b^2,(b*c)^2,
c^2,d^2,e*d^-1*c*e>;
> lS:=LowIndexSubgroups(S,<4,4>);
> for j in [1 .. #lS] do
for> print AQInvariants(lS[j]);
for> end for;
[ 0 ]
[ 2, 26 ]
[ 2, 4 ]
[ 2, 0 ]
[ 2, 2 ]
[ 2, 2 ]
[ 2, 2, 2 ]
\end{verbatim}

\bibliography{refs}

\bigbreak\noindent Department of Mathematics,
\smallbreak\noindent Rice University,
\smallbreak\noindent Houston, TX, 77005,
\smallbreak\noindent Email: tcw@rice.edu

\end{document}